\documentclass{amsart}
\newtheorem{theorem}{Theorem}[section]

\newtheorem{proposition}[theorem]{Proposition}
\newtheorem{corollary}[theorem]{Corollary}

\theoremstyle{definition}
\newtheorem{definition}[theorem]{Definition}

\newtheorem{notation}[theorem]{Notation}

\theoremstyle{remark}

\newcommand{\End}{\mathop{\mathrm{End}}\nolimits}
\newcommand{\Der}{\mathop{\mathrm{Der}}\nolimits}
\newcommand{\Ad}{\mathop{\mathrm{Ad}}\nolimits}
\newcommand{\ad}{\mathop{\mathrm{ad}}\nolimits}
\newcommand{\Hom}{\mathop{\mathrm{Hom}}\nolimits}
\newcommand{\Clif}{\mathop{\mathrm{Clif}}\nolimits}
\newcommand{\spin}{\mathop{\mathrm{spin}}\nolimits}
\newcommand{\grad}{\mathop{\mathrm{grad}}\nolimits}

\newcommand{\cL}{{\mathcal L}}
\newcommand{\cE}{{\mathcal E}}
\newcommand{\cS}{{\mathcal S}}
\newcommand{\cH}{{\mathcal H}}
\newcommand{\cT}{{\mathcal T}}
\newcommand{\bC}{{\mathbb C}}
\newcommand{\bR}{{\mathbb R}}
\newcommand{\fG}{{\mathfrak g}}
\newcommand{\fK}{{\mathfrak k}}
\newcommand{\fM}{{\mathfrak m}}

\numberwithin{equation}{section}

\begin{document}

\title[A global view of equivariant bundles and Dirac 
operators]{A global view of equivariant vector bundles \\
and Dirac operators \\ 
on some compact homogeneous spaces}

\author{Marc A. Rieffel}
\address{Department of Mathematics\\
University of California\\
Berkeley, CA\ \ 94720-3840}
\curraddr{}
\email{rieffel@math.berkeley.edu}
\thanks{The research reported here was supported in part by 
National Science Foundation Grant DMS-0500501.}
\dedicatory{Dedicated to the memory of George W. Mackey}


\subjclass[2000]{Primary 53C30; Secondary 46L87, 58J60, 
53C05}
\keywords{homogeneous spaces, equivariant vector bundles, 
connections, Dirac operators}

\begin{abstract}
In order to facilitate the comparison of Riemannian 
homogeneous spaces of compact Lie groups with noncommutative 
geometries (``quantizations'') that approximate them, we 
develop here the basic facts concerning equivariant vector 
bundles and Dirac operators over them in a way that uses only 
global constructions and arguments.  Our approach is quite 
algebraic, using primarily the modules of cross-sections of 
vector bundles.  We carry the development through the 
construction of Hodge--Dirac operators.  The inducing 
construction is ubiquitous.
\end{abstract}

\maketitle

\section{Introduction}
\label{sec1}

In the literature of theoretical high-energy physics one 
finds statements such as ``matrix algebras converge to the 
sphere'', and ``these vector bundles on the matrix algebras 
are the monopole bundles that correspond to the monopole 
bundles on the sphere''.  I have provided suitable 
definitions and theorems \cite{R6,R7} that give a precise 
meaning to the first of these statements (and I am developing 
stronger versions); but I have not yet provided precise 
meaning to the second statement, though I laid much 
groundwork for doing this in \cite{R17}.  To quantitatively 
compare vector bundles on an ordinary space, such as the 
sphere, with ``vector bundles'', i.e., projective modules, 
over a related noncommutative ``space'', it is technically 
desirable to have a description of the ordinary vector 
bundles that is as congenial to the methods of noncommutative 
geometry as possible.  This means, for example, avoiding any 
use of local coordinates, and working with the modules of 
continuous cross-sections of a bundle rather than with the 
points of the bundle itself.

The purpose of this paper is to give such a congenial 
approach for the case of equivariant vector bundles over 
homogeneous spaces of compact connected Lie groups.  (This 
includes the monopole bundles over the $2$-sphere.)  
So this paper can be viewed as largely expository, 
with its novelty being primarily in our 
presentation of the known results, and the relative 
simplicity that our approach brings to this topic.  (Compare 
with \cite{Nmz, KbN, Slb, BIL}.)  We are able to work 
entirely just in terms of functions on the Lie group.  The 
only differential geometry that we need is that which 
involves how elements of the Lie algebra give vector fields 
on the group that can differentiate functions.  Beyond this, 
our treatment is quite self-contained, and very algebraic in 
nature (somewhat along the lines found in Chapter~1 of 
\cite{Hlg}  --  see the comment near the top of page 86 of
\cite{Hlg} as to why this kind of approach is appropriate).  
The main new result is Theorem \ref{th8.4} in which we give
necessary and sufficient conditions for the Hodge-Dirac operator
corresponding to a not-necessarily torsion-free connection
to be formally self-adjoint.
Notable is our use of ``standard module 
frames'' in various places where traditionally one would have 
used local arguments involving coordinate charts.  Various 
steps of our development can be generalized in various 
directions, but for simplicity of exposition we do not 
explore these generalizations.  I expect to use the material 
presented here for my study of the quantitative relations 
between ordinary and ``quantum'' vector bundles, along the 
line that can be inferred from \cite{R17}.

In Section~\ref{sec2} of this paper we introduce equivariant 
vector bundles over homogeneous spaces and their modules of 
continuous cross-sections, while in Section~\ref{sec3} we 
give a concrete description of the algebra of module 
endomorphisms of a cross-section module.  Section~\ref{sec4} 
is devoted to discussing the tangent bundle of a homogeneous 
space through its module of cross-sections.  
Sections~\ref{sec2} to \ref{sec4} form in part a more 
leisurely version of section~$13$ and proposition~$14.3$ of 
\cite{R17}.  In Section~\ref{sec5} we discuss connections on 
equivariant vector bundles, while in Section~\ref{sec6} we 
discuss the Levi--Civita connection for an invariant 
Riemannian metric on a homogeneous space.  Finally, 
Section~\ref{sec7} is devoted to the Clifford bundle over the 
tangent space of a homogeneous space with Riemannian metric, 
in preparation for the discussion of the Hodge--Dirac 
operator that we give in Section~\ref{sec8}. One notable
aspect of Section~\ref{sec8} is that we give simple examples of
non-torsion-free connections whose Hodge-Dirac operators
are nevertheless formally self-adjoint. I have not seen this
possibility discussed in the literature.

I have tried to put in enough detail so that this paper will 
be accessible to those who have not previously met vector 
bundles and Dirac operators.

Equivariant vector bundles are very closely related to the 
induced representations that were so central to much of the 
research of George Mackey.  The inducing construction appears 
everywhere in this paper.  As an undergraduate I had the 
pleasure of taking a year-long course on projective geometry 
taught by George Mackey, and a decade later his research on 
induced representations became of great importance for my own 
research.

I developed part of the material presented here during part 
of a ten-week visit at the Isaac Newton Institute in 
Cambridge, England, in the Fall of 2006.  I am very 
appreciative of the quite stimulating and enjoyable 
conditions provided by the Newton Institute.


\section{Induced vector bundles}
\label{sec2}

In this section we assume that $G$ is a compact group, and 
that $K$ is a closed subgroup of $G$.  Then $G$ acts on the 
coset space $G/K$, which has its natural compact quotient 
topology from $G$.  We let $A = C(G/K)$, the $C^*$-algebra of 
continuous functions on $G/K$ with pointwise operations and 
supremum norm $\|\cdot\|_{\infty}$.  We will not specify 
whether the functions have values in $\bR$ or $\bC$, as 
either it will be possible to infer which from the context, 
or it does not matter for what is being discussed.  We will often
view $A$ as the subalgebra of $C(G)$ consisting of functions 
$f$ that satisfy $f(xs) = f(x)$ for $x \in G$, $s \in K$.  We 
will in general let $\lambda$ denote the action of $G$ by
left-translation on various types of functions on $G$ (or $G/K$).
In particular, we
let $\lambda$ denote the action of $G$ on $A$, defined by $(\lambda_yf)(x) = f(y^{-1}x)$ for 
$f \in A$, $y,x \in G$.

When $G$ is a Lie group, so that $G/K$ is a smooth manifold 
\cite{Wrn}, everything in this section has an evident version for 
smooth functions, but we will not state these versions, 
although we will use them in later sections.

Let $(\pi,\cH)$ be a finite-dimensional representation of 
$K$.  Since $K$ is compact we can equip the vector space 
$\cH$ with a $\pi$-invariant inner product.  Thus we will 
assume throughout that $(\pi,\cH)$ is an orthogonal or 
unitary representation.  We set
\[
\Xi_{\pi} = \{\xi \in C(G,\cH): \xi(xs) = \pi_s^{-1}(\xi(x)) 
\mbox{ for } x \in G,\ s \in K\}.
\]
Here $C(G,\cH)$ denotes the vector space of continuous 
functions from $G$ to $\cH$.  It is easily checked that 
$\Xi_{\pi}$ is in fact a module over $A$ for pointwise 
operations.  In non-commutative geometry it is often most 
convenient to put module actions on the right, so that 
operators can be put on the left.  We will follow this 
practice here for ordinary vector bundles, as done also in 
\cite{GVF}, since it is convenient here also.  Thus we view 
$\Xi_{\pi}$ as a right $A$-module.  For consistency we must 
then write scalars on the right of elements of $\cH$, and we 
have $(\xi f)(x) = \xi(x)f(x)$ for $\xi \in \Xi_{\pi}$, $f 
\in A$, $x \in G$.

The left action of $G$ on itself gives an action of $G$ on 
$C(G,\cH)$, and it is easily verified that this action 
carries $\Xi_{\pi}$ into itself.  We denote this action again 
by $\lambda$, so that $(\lambda_y\xi)(x) = \xi(y^{-1}x)$.  
Then we have the ``covariance relation'' $\lambda_y(\xi f) = 
(\lambda_y\xi)(\lambda_yf)$.  With abuse of terminology we 
set:

\begin{definition}
\label{def2.1}
The $A$-module $\Xi_{\pi}$ with its $G$-action is called the 
{\em equivariant vector bundle over $G/K$ induced from 
$(\pi,\cH)$}.
\end{definition}

We observe that for $\xi \in \Xi_{\pi}$ and $s \in K$ we have 
$(\lambda_s\xi)(e) = \xi(s^{-1}) = \pi_s(\xi(e))$, where $e$ 
is the identity element of $G$.  Thus as long as we remember how $\Xi_{\pi}$ is a space of functions on $G$, we can in this way  recover 
the original representation $(\pi,\cH)$ from $\Xi_{\pi}$ with 
its $G$-action.

The inner product on $\cH$ determines a canonical bundle 
metric (often called a Riemannian or Hermitian metric) on 
$\Xi_{\pi}$, that is, an $A$-valued inner product \cite{Lnc}, 
defined by
\[
\langle\xi,\eta\rangle_A(x) = 
\langle\xi(x),\eta(x)\rangle_{\cH}.
\]
We take our inner product on $\cH$ to be linear in its second 
variable.  It is easy to check that the $A$-valued inner 
product on $\Xi_{\pi}$ is $G$-invariant in the sense that
\[
\lambda_y(\langle\xi,\eta\rangle_A) = 
\langle\lambda_y\xi,\lambda_y\eta\rangle_A.
\]

On $G/K$ there is a unique $G$-invariant probability measure. 
 Integrating functions on $G/K$ against this measure is the 
same as viewing the functions as defined on $G$ and 
integrating them against the Haar measure on $G$ that gives 
$G$ unit mass.  Throughout this paper whenever we integrate 
over $G/K$ or $G$ it is with respect to these normalized 
measures.

On $\Xi_{\pi}$ we can define an ordinary inner product, 
$\langle\cdot,\cdot\rangle$, by
\[
\langle\xi,\eta\rangle = \int_{G/K} \langle\xi,\eta\rangle_A.
\]
The action $\lambda$ of $G$ on $\Xi_{\pi}$ preserves this 
inner product.  The completion of $\Xi_{\pi}$ for this inner 
product is the Hilbert space for Mackey's induced 
representation of $G$ from the representation $(\pi,\cH)$ of 
$K$, with the representation of $G$ just being the extension 
of $\lambda$ to the completion.  The pointwise action of $A$ 
on $\Xi_{\pi}$ extends to an action on the completion, and 
this action can be viewed as Mackey's ``system of 
imprimativity'' for the induced representation.

It can be shown that $\Xi_{\pi}$ is the space of continuous 
cross-sections of the usual equivariant vector bundle induced 
from $(\pi,\cH)$.  But we do not need this fact, though our 
development can be useful in proving this fact.  A theorem of 
Swan \cite{Swn, GVF} says that the module of continuous 
cross-sections will always be a projective module (finitely 
generated, but in this paper whenever we say ``projective'' 
we always also mean ``and finitely generated''), and 
conversely.  We include here a direct proof, taken from 
\cite{R17} but with antecedents in \cite{Lnd1}, that 
$\Xi_{\pi}$ is a projective module.  One important feature of 
the proof is that, as we will explain, it provides a way of 
obtaining a projection $p$ in a matrix algebra $M_n(A)$ for 
some $n$, such that $\Xi_{\pi}$ is isomorphic to the right 
$A$-module $pA^n$.  Such projections are crucial for 
quantifying the relation between vector bundles on compact 
metric spaces that are close together for Gromov--Hausdorff 
distance, as seen in \cite{R17}.

\setcounter{theorem}{1}
\begin{proposition}
\label{prop2.2}
For $G$, $K$ and $(\pi,\cH)$ as above, the induced module 
$\Xi_{\pi}$ is a projective $A$-module.
\end{proposition}

\begin{proof}
Find a finite-dimensional orthogonal or unitary 
representation $({\tilde \pi},{\tilde \cH})$ of $G$ such that 
$\cH$ is a subspace of ${\tilde \cH}$ and the restriction of 
${\tilde \pi}$ to $K$, acting on $\cH$, is $\pi$.  Such 
representations $({\tilde \pi},{\tilde \cH})$ exist according 
to the Frobenius reciprocity theorem, which has an elementary 
proof in our context \cite{BrD}.  Note that $C(G/K,{\tilde 
\cH})$ is a free $A$-module, with basis coming from any basis 
for ${\tilde \cH}$.  Define $\Phi$ from $\Xi_{\pi}$ to 
$C(G/K,{\tilde \cH})$ by
\[
(\Phi\xi)(x) = {\tilde \pi}_x(\xi(x))
\]
for $x \in G$.  Clearly $\Phi$ is an injective $A$-module 
homomorphism.  Let $P$ be the orthogonal projection from 
${\tilde \cH}$ onto $\cH$, and define a function $p$ from $G$ 
to $\cL({\tilde \cH})$, the algebra of linear operators on 
${\tilde \cH}$, by
\[
p(x) = {\tilde \pi}_xP {\tilde \pi}_x^*.
\]
It is easily seen that $p$ is a projection in the algebra 
$C(G/K,\cL({\tilde \cH}))$.  Furthermore, this algebra acts 
as endomorphisms on the free $A$-module $C(G/K,{\tilde \cH})$ 
in the evident pointwise way, and it is easy to check 
\cite{R17} that $p$ is the projection onto the range of 
$\Phi$.  Thus the range of $\Phi$, and so $\Xi_{\pi}$, is 
projective.
\end{proof}

We remark that the actual vector bundle corresponding to 
$\Xi_{\pi}$ can be viewed as assigning to each point ${\dot 
x}$ of $G/K$ the range subspace of $p(x)$.  Notice that for a 
given $(\pi, \cH)$ there may be many choices for the 
representation $({\tilde \pi},{\tilde \cH})$ above, and so 
many choices of the projection $p$.

In case $G$ is a Lie group, it is known that 
finite-dimensional representations (as homomorphisms between 
Lie groups) are real-analytic, so smooth.  Consequently the 
projection $p$ of the above proof is smooth, and this shows 
that the subspace $\Xi_{\pi}^{\infty}$ of smooth elements of 
$\Xi_{\pi}$ is a projective module over $C^{\infty}(G/K)$.

Let $A$ be any unital $C^*$-algebra, and let $\Xi$ be a right 
$A$-module.  Assume that $\Xi$ is equipped with an $A$-valued 
inner product \cite{Lnc}.  A finite sequence $\{\eta_j\}$ of 
elements of $\Xi$ is said to be a {\em standard module frame} 
for $\Xi$ if the ``reproducing formula''
\[
\xi = \sum \eta_j\langle\eta_j,\xi\rangle_A
\]
holds for all $\xi \in \Xi$.  (It is easily seen, p.~$46$ of 
\cite{R17}, that if $\Xi$ has a standard module frame, 
$\{\eta_j\}_{j=1}^n$, then $\Xi$ is a projective $A$-module, 
and the matrix $p = \{\langle\eta_j,\eta_k\rangle_A\}$ in 
$M_n(A)$ is a projection such that $\Xi \cong pA^n$.)  We can 
construct standard module frames for the induced vector 
bundles $\Xi_{\pi}$ as follows.  For $({\tilde \pi},{\tilde 
\cH})$ as used above, let $\{e_j\}$ be an orthonormal basis 
for ${\tilde \cH}$.  
For each $j$ define $\eta_j$ by $\eta_j(x) = P{\tilde 
\pi}_x^{-1}e_j$ for $x \in G$.  Each $\eta_j$ is easily seen 
to be in $\Xi_{\pi}$.  Then for $\xi \in \Xi_{\pi}$ we have
\begin{eqnarray*}
(\sum \eta_j(x)\langle\eta_j,\xi\rangle_A)(x) &= &\sum 
P{\tilde \pi}_x^{-1}e_j\langle P{\tilde 
\pi}_x^{-1}e_j,\xi(x)\rangle_{\cH} \\
&= & P\sum ({\tilde \pi}_x^{-1}e_j) \langle ({\tilde 
\pi}_x^{-1}e_j),\xi(x)\rangle_{\cH} = P\xi(x) = \xi(x),
\end{eqnarray*}
where we have used that $\{{\tilde \pi}_x^{-1}e_j\}$ is 
equally well an orthonormal basis for ${\tilde \cH}$.  Thus 
$\{\eta_j\}$ has the reproducing property, and so is a 
standard module frame for $\Xi_{\pi}$.  (We remark that it 
may happen that $\eta_j = 0$ for certain $j$'s.)  We will 
make good use of standard module frames in some of the next 
sections.


\section{Endomorphism bundles}
\label{sec3}

In this section we will give a description of the 
endomorphism bundles of induced bundles.  This description 
will be useful in our later discussion of connections.  We 
will continue to work just with continuous functions and 
cross-sections, and we will leave it to the reader to notice 
that when $G$ is a Lie group everything said in this section 
has a smooth version.  We will need these smooth versions in 
later sections.

As before, let $(\pi,\cH)$ be a finite-dimensional orthogonal 
or unitary representation of $K$, and let $\Xi_{\pi}$ be the 
induced bundle.  We want to describe $\End_A(\Xi_{\pi})$ in 
terms of $(\pi,\cH)$.  Let $\{\eta_j\}$ be a standard module 
frame for $\Xi_{\pi}$.  For any $T \in \End_A(\Xi_{\pi})$ and 
any $\xi \in \Xi_{\pi}$ we have
\begin{eqnarray*}
T\xi &= &T\left(\sum \eta_j\langle\eta_j,\xi\rangle_A\right) 
= \sum (T\eta_j)\langle\eta_j,\xi\rangle_A \\
&= &\left(\sum \langle T\eta_j,\eta_j\rangle_E\right)\xi,
\end{eqnarray*}
where for $\zeta,\eta \in \Xi_{\pi}$ we let 
$\langle\zeta,\eta\rangle_E$ denote the ``rank-one'' operator 
on $\Xi_{\pi}$ defined by $\langle\zeta,\eta\rangle_E\xi = 
\zeta\langle\eta,\xi\rangle_A$.  Then 
$\langle\cdot,\cdot\rangle_E$ is an inner product on 
$\Xi_{\pi}$ with values in $E = \End_A(\Xi_{\pi})$.  It is 
linear in the first variable, and elements of $E$ pull out of 
the first variable as if they were scalars \cite{Lnc}; we 
view $\Xi_{\pi}$ as a left $E$-module. The calculation above 
shows that for any $T \in E$ we have $T = \sum \langle 
T\eta_j,\eta_j\rangle_E$.  Thus $E$ is spanned by the 
``rank-one'' operators.

For any $\zeta,\eta,\xi \in \Xi_{\pi}$ we have
\[
(\langle\zeta,\eta\rangle_E\xi)(x) = 
\zeta(x)\langle\eta(x),\xi(x)\rangle_{\cH} = 
\langle\zeta(x),\eta(x)\rangle_0(\xi(x)),
\]
where now $\langle\cdot,\cdot\rangle_0$ denotes the usual 
ordinary rank-one operator on $\cH$ given by two vectors.  
Since we just saw that the ``rank-one'' operators 
$\langle\zeta,\eta\rangle_E$ span $E$, we see now that every 
operator in $E$ is given by a function in $C(G,\cL(\cH))$, as 
we would expect.  Furthermore, for any $x \in G$ and $s \in 
K$ we have
\[
\langle\zeta(xs),\eta(xs)\rangle_0 = 
\langle\pi_s^{-1}(\zeta(x)),\pi_s^{-1}(\eta(x))\rangle_0 = 
\pi_s^{-1}\langle\zeta(x),\eta(x)\rangle_0\pi_s.
\]
It follows that if $T \in C(G,\cL(\cH))$ represents an 
element of $E$ then $T(xs) = \pi_s^{-1}T(x)\pi_s$.  But a 
simple direct check shows that any $T$ satisfying this 
property gives an element of $E$.  Thus we have obtained

\begin{proposition}
\label{prop3.1}
We can identify $\End_A(\Xi_{\pi})$ with
\[
\cE_{\pi} = \{T \in C(G,\cL(\cH)): T(xs) = 
\pi_s^{-1}T(x)\pi_s \mbox{ for } x \in G,\ s \in K\},
\]
where the action of such a \ $T$ on $\xi \in \Xi_{\pi}$ is 
given just by $(T\xi)(x) = T(x)(\xi(x))$.
\end{proposition}

There is an evident action, say $\alpha$, of $K$ on 
$\cL(\cH)$, given by $\alpha_s(\tau) = \pi_s\tau\pi_s^{-1}$, 
and we see that $\cE_{\pi}$ is just the ``induced algebra over 
$G/K$ for the action $\alpha$''.  As such, $\cE_{\pi}$ is an 
$A$-module.  But it is also an $A$-module just from the fact 
that $A$ is commutative so that $A$ is contained in the 
center of $\cE_{\pi}$, and we notice that these two 
$A$-module structures coincide.


\section{The tangent bundle}
\label{sec4}

We now assume that $G$ is a compact connected Lie group with 
Lie algebra $\fG$.  We assume again that $K$ is a closed 
subgroup of $G$, which need not be connected, and we denote 
its Lie algebra (for its connected component of the identity 
element $e$) by $\fK$.  Then $G/K$ is a compact smooth 
manifold \cite{Wrn}.  We seek an attractive realization of 
its tangent bundle, as a projective $A$-module, where now we 
set $A = C^{\infty}(G/K)$, with functions being real-valued.  
We use below some standard facts \cite{Hlg, Wrn} about Lie 
groups and their Lie algebras.  This section and the next two 
have a number of points of contact with section~$5$ of 
\cite{BIL}.

Much as before, we view $C^{\infty}(G/K)$ as a subalgebra of 
$C^{\infty}(G)$.  The action $\lambda$ of $G$ on 
$C^{\infty}(G)$ and $C^{\infty}(G/K)$ gives an infinitesimal 
action of $\fG$, defined for $X \in \fG$ by $(\lambda_Xf)(x) 
= (Xf)(x) = D_0^t(f(\exp(-tX)x))$, where we write $D_0^t$ for 
$(d/dt)|_{t=0}$.  This gives a Lie algebra homomorphism from 
$\fG$ into the Lie algebra $\Der(C^{\infty}(G))$ of 
derivations of the algebra $C^{\infty}(G)$, and so into 
$\Der(C^{\infty}(G/K))$. (See proposition 2.4 of chapter 0 of 
\cite{Tyl}.) The derivations $\lambda_X$ of $C^{\infty}(G/K)$ 
are often called the {\em fundamental} vector fields for 
$C^{\infty}(G/K)$.  Since $C^{\infty}(G)$ is commutative, 
$\Der(C^{\infty}(G))$ is a module over $C^{\infty}(G)$, and a 
basis for $\fG$ gives a module basis for 
$\Der(C^{\infty}(G))$, that is, $\Der(C^{\infty}(G))$ is a 
free $C^{\infty}(G)$-module.  We can thus realize the 
elements of $\Der(C^{\infty}(G))$, which is the module of 
smooth cross-sections of the tangent bundle of $G$, as 
elements of $C^{\infty}(G,\fG)$, where the action of $W \in 
C^{\infty}(G,\fG)$ on $C^{\infty}(G)$ is given by
\[
(\lambda_Wf)(x) = D_0^t(f(\exp(-tW(x))x)).
\]
Accordingly, we denote $C^{\infty}(G,\fG)$ by $\cT(G)$, for 
``tangent space''.

Now for $f \in A = C^{\infty}(G/K)$ and $W \in \cT(G)$ we 
need to have $W(xs) = W(x)$ for $x \in G$ and $s \in K$ if we 
want $\lambda_Wf \in A$.  Furthermore, if for some $x_1 \in 
G$ we have $W(x_1) \in \Ad_{x_1}(\fK)$, then 
$\Ad_{x_1}^{-1}(W(x_1)) \in \fK$, so that
\begin{eqnarray*}
(\lambda_Wf)(x_1) &= &D_0^t(f(\exp(-tW(x_1))x_1) \\
&= &D_0^t(f(x_1\exp(-t\Ad_{x_1}^{-1}(W(x_1))) = D_0^t(f(x_1)) 
= 0.
\end{eqnarray*}
Conversely, if $(\lambda_Wf)(x_1) = 0$ for all $f \in A$, 
then $W(x_1) \in \Ad_{x_1}(\fK)$.

Choose an $\Ad$-invariant inner product, 
$\langle \cdot, \ \cdot\rangle_\fG$, on $\fG$, which we 
fix for the rest of this paper.  Let $\fM$ denote the 
orthogonal complement to $\fK$ in $\fG$.  Let $P$ be the 
orthogonal projection of $\fG$ onto $\fM$.  For any $W \in 
\cT(G)$ define ${\tilde W}$ by ${\tilde W}(x) = 
\Ad_x(P(\Ad_x^{-1}(W(x))))$.  Then $W(x) - {\tilde W}(x) \in 
\Ad_x(\fK)$ for all $x$, and so $W-{\tilde W}$ acts on $A$ as 
the $0$-derivation.  This and the earlier calculations 
suggest that we consider
\begin{equation}
\label{eq4.1}
\{W \in C^{\infty}(G,\fG): W(x) \in \Ad_x(\fM) \mbox{ and } 
W(xs) = W(x) \mbox{ for all } x \in G,\ s \in K\}.
\end{equation}
And indeed it is easily seen that any such $W$, acting as 
derivations of $C^{\infty}(G)$ by the earlier formula, 
carries $A$ into itself, and that every smooth vector field 
on $G/K$ is represented in this way.  However, it is 
inconvenient that the range space of the $W$'s is not 
constant on $G$.  But
\[
\exp(W(x))x = x \exp(\Ad_x^{-1}(W(x))),
\]
and if $W$ is as in 4.1 then $\Ad_x^{-1}(W(x)) \in \fM$ for 
all $x$.  When we also take into account the effect on the 
right $K$-invariance, we are led to:

\setcounter{theorem}{1}
\begin{notation}
\label{note4.2}
With notation as above, set
\[
\cT(G/K) = \{W \in C^{\infty}(G,\fM): W(xs) = 
\Ad_{s^{-1}}(W(x)) \mbox{ for } x \in G,\ s \in K\}.
\]
\end{notation}

We let elements of $\cT(G/K)$ act as derivations of 
$C^{\infty}(G/K)$ by
\[
(\delta_Wf)(x) = D_0^t(f(x\exp(tW(x)))),
\]
where we have also left behind the earlier minus sign.  It is 
clear that $\cT(G/K)$ is a module over $A$ for pointwise 
operations.  We recognize $\cT(G/K)$ as just the induced 
bundle for the representation $\Ad$ restricted to $K$ on 
$\fM$.  It is easy to check that $\cT(G/K)$ does give 
derivatives in all tangent directions at any point of $G/K$. We
can then to use the fact that $\cT(G/K)$ is an $A$-module to 
show that it does contain all the smooth vector fields on 
$G/K$. Thus it does represent the space of smooth 
cross-sections of the tangent bundle of $G/K$, though we will 
not explicitly need this fact.  This description of the 
tangent bundle can be found, for example, in \cite{Slb, Frd, 
BIL}.

For $X \in \fG$ the corresponding fundamental vector field on 
$G/K$ is given, in the form used for \eqref{eq4.1}, by
\[
X(x) = \Ad_x(P \Ad_x^{-1}(X)).
\]
Then in the form used for the definition of $\cT(G/K)$ this 
fundamental vector field, which we now denote by ${\hat X}$, 
is given by
\[
{\hat X}(x) = -P \Ad_x^{-1}(X).
\]
(See section~$0.3$ of \cite{Slb}.)  Since the map from 
elements of $\fG$ to vector fields is a Lie algebra 
homomorphism, we have
\[
[{\hat X},{\hat Y}](x) = -P \Ad_x^{-1}([X,Y]),
\]
where the brackets on the left denote the commutator of 
${\hat X}$ and ${\hat Y}$ as operators on $A$.  Also, a quick 
calculation shows that $\lambda_y{\hat X} = 
(\Ad_yX)^{\wedge}$ for $y \in G$.  We remark that ${\hat X}$ 
may well take value $0$ at some points of $G/K$ --- we are 
confronting the fact that the tangent bundle of $G/K$ may 
well not be a trivial bundle.  It may even happen that for 
certain $X$'s in $\fG$ we have ${\hat X} \equiv 0$. We warn 
the reader that for general $V, W \in \cT(G/K)$ the commutator
$[V, \ W]$ as derivations of $A$ is again an element of $\cT(G/K)$ but
it is usually \emph{not} given by any pointwise formula. The torsion-free
condition shows how to express this commutator in terms of the
Levi-Civita connection or other torsion-free connections as we will
see later.

The restriction to $\fM$ of our chosen inner product on $\fG$
determines a chosen $G$-invariant Riemannian metric on $G/K$,
specializing what was done in Section \ref{sec2} .
If $\{X_j\}$ is an orthonormal basis for $\fG$, then one can 
chase through the discussion at the end of the previous 
section to see that $\{{\hat X}_j\}$ is a standard module 
frame for $\cT(G/K)$ for this Riemannian metric.  But this is also easy to verify 
directly:  For any $W \in \cT(G/K)$ and $x \in G$ we have
{\setlength\arraycolsep{2pt}
\begin{eqnarray*}
\left( \sum {\hat X}_j\langle{\hat X}_j,W\rangle_A\right)(x) 
&= &\sum - (P \Ad_x^{-1}X_j)\langle -P \Ad_x^{-1} X_j, \ 
W(x)\rangle_{\fG} \\
&= &P \sum (\Ad_x^{-1}X_j)\langle \Ad_x^{-1}X_j, \ 
W(x)\rangle_{\fG} = PW(x) = W(x),
\end{eqnarray*}   
}  
since $\{\Ad_x^{-1}X_j\}$ is equally well an orthonormal 
basis for $\fG$.  Of course, again we may have ${\hat X}_j 
\equiv 0$ for some $j$'s.


\section{Connections}
\label{sec5}

In this section we will give a description of all 
$G$-invariant connections on an induced bundle.  By a 
connection (or covariant derivative) on the space $\Xi$ of 
smooth cross-sections of a vector bundle on a manifold $M$ we 
mean \cite{Hlg, GVF} a linear transformation $\nabla$ from 
the space $\cT(M)$ of smooth tangent-vector fields on $M$ to 
linear transformations on $\Xi$ such that for  $W \in 
\cT(M)$,  $\xi \in \Xi$ and $f \in C^{\infty}(M)$ we have the 
Leibniz rule
\[
\nabla_W(\xi f) = (\nabla_W\xi)f + \xi(\delta_Wf)
\]
and the rule
\[
\nabla_{Wf}(\xi) = (\nabla_W\xi)f.
\]
For the case we have been considering, in which $M = G/K$ and 
$\Xi$ is an induced vector bundle $\Xi_{\pi}$, there is a 
canonical connection, $\nabla^0$, defined by
\[
(\nabla_W^0\xi)(x) = D_0^t(\xi(x\exp(tW(x)))).
\]
To see that $\nabla_W^0\xi$ is indeed in $\Xi_{\pi}$ we 
calculate that for $x \in G$ and $s \in K$
\begin{eqnarray*}
(\nabla_W^0\xi)(xs) &= &D_0^t(\xi(xs\exp(tW(xs)))) 
= D_0^t(\xi(xs\exp(t \Ad_s^{-1}(W(x))))) \\
&= &D_0^t(\xi(x\exp(t(W(x)))s)) = 
\pi_s^{-1}((\nabla_W^0\xi)(x)),
\end{eqnarray*}
as needed.  The other properties stated above for the 
definition of a connection are easily verified by similar 
calculations.

There is a standard definition \cite{Hlg} of what it means 
for a connection $\nabla$ on an equivariant vector bundle to 
be invariant for the group action, and when this definition 
is applied to our induced bundles $\Xi_{\pi}$ over $G/K$, it 
requires that
\[
\lambda_y(\nabla_W\xi) = \nabla_{\lambda_yW}(\lambda_y\xi)
\]
for all $y \in G$, $W \in \cT(G/K)$ and $\xi \in \Xi_{\pi}$.  
It is straightforward to check that the canonical connection 
$\nabla^0$ is $G$-invariant.

Let $\nabla$ be any connection on $\Xi_{\pi}$, and set $L = 
\nabla - \nabla^0$.  It is easily checked that $L_W \in 
\End_A(\Xi_{\pi})$ for each $W \in \cT(G/K)$, and so $L_W$ 
can be represented as a function in $\cE_{\pi}$ according to 
Proposition~\ref{prop3.1}.  We mentioned in 
Section~\ref{sec3} that $\cE_{\pi}$ is an $A$-module because 
$A$ is commutative.  It is easily checked that $L$ is an 
$A$-module homomorphism from $\cT(G/K)$ to $\cE_{\pi}$, and 
that, conversely, if $L$ is any $A$-module homomorphism from 
$\cT(G/K)$ to $\cE_{\pi}$ then $\nabla^0 + L$ is a connection 
on $\Xi_{\pi}$.  This is just the well-known fact that the 
set of connections forms an affine space over 
$\Hom_A(\cT(G/K),\End_A(\Xi_{\pi}))$.
If $\nabla$ is $G$-invariant, then, because $\nabla^0$ also 
is $G$-invariant, so is $L$, in the sense that 
$\lambda_y(L_W\xi) = L_{\lambda_yW}(\lambda_y\xi)$ for all $y 
\in G$, $W \in \cT(G/K)$ and $\xi \in \Xi_{\pi}$.  When we 
view $L_W$ as a function in $\cE_{\pi}$, invariance gives
\[
L_W(y^{-1}x)\xi(y^{-1}x) = (\lambda_y(L_W\xi))(x) = 
(L_{\lambda_yW}(\lambda_y\xi))(x) = 
L_{\lambda_yW}(x)\xi(y^{-1}x).
\]
Since this is true for all $\xi$, we see that for all $y \in 
G$ we have
\[
\lambda_yL_W = L_{\lambda_yW}
\]
as functions on $G$.  In particular $L_W(y^{-1}) = 
(\lambda_yL_W)(e) = L_{\lambda_yW}(e)$, so that $L$ is 
determined once we know $L_W(e)$ for all $W \in \cT(G/K)$.  
Note further that since $L_W \in \cE_{\pi}$, we have by 
invariance
\[
L_{\lambda_sW}(e) = L_W(es^{-1}) = \pi_sL_W(e)\pi_s^{-1}
\]
for $s \in K$.

Suppose now that $\{X_j\}$ is an orthonormal basis for $\fG$, 
so that $\{{\hat X}_j\}$ is a standard module frame for 
$\cT(G/K)$, as seen near the end of Section~\ref{sec4}.  Then
\[
L_W(e) = \sum L_{{\hat X}_j}(e)\langle{\hat 
X}_j,W\rangle_A(e),
\]
so that $L_W(e)$ is determined once we know $L_{{\hat 
X}_j}(e)$ for all $j$. But ${\hat X}(e) = 0$ if $X \in \fK$.
Thus in the above expression we only need to sum over a
basis for $\fM$. Equivalently, $L_W(e)$ is determined once we know
$L_{\hat X}(e)$ for all $X \in \fM$. Notice that on $\fM$ the map
$X \mapsto {\hat X}$ is injective. Define $\gamma_L$ on $\fM$ by
$\gamma_L(X) = L_{\hat X}(e)$, so that $L$ is determined by
$\gamma_L$.
Then $\gamma_L$ is a real-linear transformation from 
$\fM$ to $\cL(\cH)$.  
Recall that $\lambda_y({\hat 
X}) = (\Ad_y(X)) \ \hat{} \ $.  Then from the property of $L_W(e)$ 
obtained at the end of the previous paragraph we have
\[
\gamma_L(\Ad_s(X)) = \pi_s\gamma_L(X)\pi_s^{-1}
\]
for $s \in K$ and $X \in \fM$.
The following 
theorem has its roots at least back in Nomizu \cite{Nmz}.  
See also section~X.2 of \cite{KbN} and section~$0$ of 
\cite{Slb}.

\begin{theorem}
\label{th5.1}
With notation as above, the map $L \mapsto \gamma_L$ gives a 
bijection between the set of $G$-invariant connections on 
$\Xi_{\pi}$ and the set of linear operators $\gamma$ from 
$\fM$ to $\cL(\cH)$ with the property that
\[
\gamma(\Ad_s(X)) = \pi_s\gamma(X)\pi_s^{-1}
\]
for all $s \in K$ and $X \in \fM$.
\end{theorem}

\begin{proof}
We have shown above that every invariant connection $\nabla$ 
gives rise to $L = \nabla - \nabla^0$, and $L$ gives rise to 
$\gamma_L$, which in turn determines $L$ and so $\nabla$.  We 
must show, conversely, that any $\gamma$ as in the statement 
of the theorem gives rise to an $L \in 
\Hom_A(\cT(G/K),\cE_{\pi})$ such that $\gamma = \gamma_L$, 
and so gives rise to the connection $\nabla = \nabla^0 + L$, 
which is $G$-invariant.  We first define $L$ on each $\hat X$ 
for $X \in \fG$ by
\[
L_{\hat X}(x) = \gamma(P(\Ad_x^{-1}(X))).
\]
It is easy to check that $L_{\hat X} \in \cE_{\pi}$, that is, that
\[
L_{\hat X}(xs) = \pi_s^{-1}L_{\hat X}(x)\pi_s
\]
for $x \in G$ and $s \in K$.  We then choose a standard 
module frame $\{{\hat X}_j\}$ for $\cT(G/K)$, and set $L_W = 
\sum L_{{\hat X}_j}\langle{\hat X}_j,W\rangle_A$.  It is then 
easy to check that this $L$ has the desired properties,
and that $\gamma_L = \gamma$.
\end{proof}

In the presence of a bundle metric on a vector bundle, a 
connection is said to be compatible with the bundle metric if 
the Leibniz rule
\[
\delta_W(\langle\xi,\eta\rangle_A) = 
\langle\nabla_W\xi,\eta\rangle_A + 
\langle\xi,\nabla_W \eta\rangle_A
\]
holds.  It is easy to verify that the canonical connection 
$\nabla^0$ on $\Xi_{\pi}$ is compatible with the bundle 
metric that we have been using.  If $\nabla$ is another 
connection on $\Xi_{\pi}$ that is compatible, and if we set 
$L = \nabla - \nabla^0$, then we see that $L$ must satisfy
\[
0 = \langle L_W\xi,\eta\rangle_A + 
\langle\xi,L_W\eta\rangle_A.
\]
This says that $L$, as a function with values in $\cL(\cH)$, 
must in fact have its values in the subspace $\cL^{sk}(\cH)$ 
of skew-symmetric or skew-Hermitian operators on $\cH$.  When 
$L$ is $G$-invariant and we define $\gamma_L$ as above by 
$\gamma_L(X) = L_{\hat X}(e)$, this then says exactly that 
$\gamma_L$ must have its values in $\cL^{sk}(\cH)$.  Thus we 
obtain:

\begin{corollary}
\label{cor5.2}
With notation as above, the map $L \mapsto \gamma_L$ gives a 
bijection between the set of $G$-invariant compatible 
connections on $\Xi_{\pi}$ and the set of linear operators 
$\gamma$ from $\fM$ to $\cL^{sk}(\cH)$ with the 
property that
\[
\gamma(\Ad_s(X)) = \pi_s\gamma(X)\pi_s^{-1}
\]
for all $s \in K$ and $X \in \fM$.
\end{corollary}


\section{The Levi--Civita connection}
\label{sec6}

For a connection $\nabla$ on a tangent bundle itself it makes 
sense to talk about its torsion \cite{Hlg, GVF}, which for a 
connection $\nabla$ on $\cT(G/K)$ is the bilinear form 
$T_{\nabla}$ defined by
\[
T_{\nabla}(V,W) = \nabla_V(W) - \nabla_W(V) - [V,W]
\]
for $V,W \in \cT(G/K)$, with values in $\cT(G/K)$.  Then the 
Levi--Civita connection associated to a Riemannian metric on 
the tangent space is by definition, for our case of 
$\cT(G/K)$, the (necessarily unique) connection $\nabla$ 
which is compatible with the Riemannian metric and has 
$T_{\nabla} \equiv 0$.  Since our Riemannian metric is 
$G$-invariant, we can expect its Levi--Civita connection to 
be $G$-invariant also.

Let us calculate the torsion of the canonical connection 
$\nabla^0$ on $\cT(G/K)$.  
Now for any connection $\nabla$ it 
is not difficult to verify that $T_{\nabla}$ is $A$-bilinear. 
 (See \S8 of chapter~$1$ of \cite{Hlg}.)  Thus for the 
reasons seen earlier, it is sufficient to calculate with 
fundamental vector fields.  
Now to begin with, for $X, Y \in \fG$ we have
\setlength\arraycolsep{2pt}
\begin{eqnarray*}
(\nabla_{\hat X}^0({\hat Y}))(x) &= &D_0^t {\hat 
Y}(x \exp(-tP\Ad^{-1}_x (X))) \\
&= &D_0^t(-P(\Ad_{\exp(tP\Ad^{-1}_x(X)} \Ad^{-1}_x(Y))) 
= -P([P\Ad^{-1}_x(X), \Ad^{-1}_x(Y)]).
\end{eqnarray*}
If a connection $\nabla$ is $G$-invariant, 
then it is easily seen that $T_{\nabla}$ is also, in the 
sense that
\[
\lambda_y(T_{\nabla}(V,W)) = 
T_{\nabla}(\lambda_yV,\lambda_yW).
\]
Then it suffices to calculate at $e$. Accordingly
\[
T_{\nabla^0}({\hat X},{\hat Y})(e) = -P([PX,Y]) + P([PY,X]) - 
[{\hat X},{\hat Y}](e).
\]
Since $X \to {\hat X}$ is a Lie algebra homomorphism, $[{\hat 
X},{\hat Y}](e) = [X,Y]^{\wedge}(e) = -P([X,Y])$, and so if 
we let $Q = I - P$,
\begin{eqnarray*}
T_{\nabla^0}({\hat X},{\hat Y})(e) &= &P([PY,X] - [PX,Y] + 
[X,Y]) 
= P([PY,X] + [QX,Y]) \\
&= &P(-[X,PY] + [QX,PY] + [QX,QY]) = -P([PX,PY]),
\end{eqnarray*}
since $[QX,QY] \in \fK$ so that $P([QX,QY]) = 0$.  
From this we find easily that for $V, W \in \cT(G/K)$ we have
\[
T_{\nabla^0}(V, W)(x) = -P[V(x), W(x)]  \ .
\]

We thus 
see that $T_{\nabla^0} \equiv 0$ exactly if $[\fM,\fM] 
\subseteq \fK$, which is exactly the condition for $G/K$ to 
be a symmetric space \cite{Hlg, Frd}, since the involutive 
transformation on $\fG$ which is the identity on $\fK$ and 
the negative of the identity on $\fM$ is then a Lie algebra 
homomorphism.  We thus find the well-known fact (see 
section~$3.5$ of \cite{Frd})  that $\nabla^0$ is the 
Levi--Civita connection exactly if $G/K$ is a symmetric 
space.

If $G/K$ is not a symmetric space, then the Levi--Civita 
connection will be 
of the form $\nabla^0 + L$.  We can again work at $e$, and 
one soon sees 
that if we define $L^0$ at $e$ by
\[
(L^0_{\hat X}({\hat Y}))(e) = (1/2)P([PX,PY])
\]
and extend $L^0$ to $G/K$ by
\[
(L^0_{\hat X}({\hat Y}))(x) = (L^0_{\lambda_x^{-1}{\hat 
X}}(\lambda_x^{-1}{\hat Y}))(e),
\]
noting that $\lambda_x^{-1}({\hat X}) = 
(\Ad_x^{-1}(X)) \ \hat{} \ $, and 
finally extend $L^0$ to $\cT(G/K)$ by using a standard module 
frame of 
fundamental vector fields as done earlier, then 
$L^0 \in \Hom_A(\cT(G/K),\End_A(\cT(G/K))$, and
{\setlength\arraycolsep{2pt}
\begin{eqnarray*}
(T_{\nabla^0+L^0}({\hat X}, {\hat Y})(e) 
&= &(T_{\nabla^0}({\hat X},{\hat Y}))(e) + 
(L^0_{\hat X}({\hat Y}))(e) - (L_{\hat Y}^0({\hat X}))(e) \\
= -P([PX,PY]) &+ &(1/2)P([PX,PY]) - (1/2)P([PY,PX]) = 0.
\end{eqnarray*}}
Thus $\nabla^0 + L^0$ has torsion $0$.  We see that the 
$\gamma$ for $L^0$ 
as in Theorem~\ref{th5.1} is defined by
\[
\gamma_X = (1/2)P \circ \ad_{X} 
\]
for all $X \in \fM$ as an operator on $\fM$.  
(Compare with theorem~X.2.10 of 
\cite{KbN} and 
lemma~0.4.3 of \cite{Slb}.)  It is easily seen that 
$\gamma_X \in \cL^{sk}(\fM)$ because the inner product on 
$\fG$ was chosen 
to be $\Ad$-invariant, so that $\ad_Z$ is a skew-symmetric 
operator on $\fG$ 
for every $Z \in \fG$.  From Corollary~\ref{cor5.2} it 
follows that 
$\nabla^0 + L^0$ is compatible with the canonical Riemannian 
metric 
on $\cT(G/K)$, so that $\nabla^0 + L^0$ is the Levi--Civita 
connection for 
that Riemannian metric.  When one carries through the 
calculations with the 
fundamental vector fields one finds that $L^0$ is given for 
general 
$V,W \in \cT(G/K)$ by
\[
(L_V^0W)(x) = (1/2)P([V(x),W(x)]).
\]
One easily checks directly that $L^0 \in 
\Hom_A(\cT(G/K),\End_A(\cT(G/K)))$ 
when $L^0$ is defined in this way, that $L^0$ is 
$G$-invariant, and is 
skew-symmetric.  But it does not seem so easy to check 
directly that 
$\nabla^0 + L^0$ has $0$ torsion, or even to guess that the 
above formula 
is the correct one for $L^0$, without working with the 
fundamental vector 
fields.  In summary: 

\begin{theorem}
\label{th6.1}
With notation as above, the Levi--Civita connection for the 
canonical  metric on $G/K$ is $\nabla^0 + L^0$ where 
$\nabla^0$ is the canonical connection on $\cT(G/K)$ and 
$L^0$ is defined by
\[
(L^0_VW)(x) = (1/2)P([V(x),W(x)]).
\]
If $[\fM,\fM] \subseteq \fK$, then $L^0 \equiv 0$.
\end{theorem}


\section{The Clifford-algebra bundle}
\label{sec7}

We can form the Clifford algebra, $\Clif(\fM)$, over $\fM$ 
with respect to the inner product on $\fM$ (coming from that 
on $\fG$).  By definition \cite{Frd, GVF} $\fM$ sits as a 
vector subspace of $\Clif(\fM)$ and generates $\Clif(\fM)$ as a 
unital algebra.  We follow the common convention in 
Riemannian geometry that the defining relation for elements 
of $\fM$ is
\[
X \cdot Y + Y \cdot X = -2\langle X,Y\rangle_{\fM}1,
\]
where we denote the product in $\Clif(\fM)$ by ``$\cdot$''.  
One can view $\Clif(\fM)$ as a deformation of the exterior 
algebra over $\fM$ in the direction of the inner product.  
Any isometric operator on $\fM$ extends uniquely to an 
algebra automorphism of $\Clif(\fM)$.  Since the action $\Ad$ 
of $K$ on $\fM$ is by isometries, it extends to an action of 
$K$ as algebra automorphisms of $\Clif(\fM)$. We denote this 
action again by $\Ad$.

The tangent space at each point of $G/K$ is isomorphic to 
$\fM$, and we can form the smooth cross-section algebra of 
the bundle of the corresponding
Clifford algebras. We denote it by 
$\Clif(\cT(G/K))$.  We seek an explicit description of this 
algebra in terms of our explicit description of $\cT(G/K)$.  
Since $\Clif(\cT(G/K))$ should be an $A$-module algebra and 
contain $\cT(G/K)$ as a submodule that generates it, 
we are led to set
\[
\Clif(\cT(G/K)) = \{\varphi \in C^{\infty}(G,\Clif(\fM)): 
\varphi(xs) = \Ad_s^{-1}(\varphi(x)) \mbox{ for } x \in G,\ s 
\in K\}.
\]
This is, of course, yet another ``induced'' algebra.  It 
contains $A = C^{\infty}(G/K)$ in its center in the evident 
way, it contains $\cT(G/K)$, and it is an algebra generated 
by $\cT(G/K)$ with the expected Clifford-algebra relations, 
namely
\[
V \cdot W + W \cdot V = -2\langle V, W \rangle_A  .
\]

Since our canonical connection is compatible with the 
Riemannian metric on $\cT(G/K)$, it extends to a connection 
on $\Clif(\cT(G/K))$, which we denote again by $\nabla^0$.  
It can be defined directly by
\[
(\nabla_W^0\varphi)(x) = D_0^t(\varphi(x\exp(tW(x))).
\]
This clearly satisfies the Leibniz rule
\[
\nabla_W^0(\varphi \cdot \psi) = (\nabla_W^0\varphi)\cdot\psi 
\ \ + \  \ \varphi\cdot(\nabla_W^0\psi)
\]
for $\varphi,\psi \in \Clif(\cT(G/K))$, that is, $\nabla_W^0$ 
is a derivation of $\Clif(\cT(G/K))$ for each $W \in 
\cT(G/K)$.  Note that when this $\nabla_W^0$ is restricted to 
$A \subset \Clif(\cT(G/K))$ it gives $\delta_W$.  Since the 
Riemannian metric on $\cT(G/K)$ is invariant for the left 
action $\lambda$ of $G$ on $\cT(G/K)$ by translation, 
$\lambda$ extends to an action of $G$ by algebra 
automorphisms on $\Clif(\cT(G))$, which we again denote by 
$\lambda$.  It is of course given by $(\lambda_y\varphi)(x) = 
\varphi(y^{-1}x)$, and when restricted to $A \subset 
\Clif(\cT(G/K))$ it gives the original action of $G$ on $A$.

The standard Dirac operator for a Riemannian manifold is 
defined in terms of the Levi--Civita connection, and so if 
$G/K$ is not a symmetric space, we should consider instead 
the connection $\nabla^0 + L^0$ defined in the previous 
section.  Recall that $L_W^0(x)$ is a skew-symmetric operator 
on $\fM$ for each $W \in \cT(G/K)$ and $x \in G$.  Now any 
skew-symmetric operator, say $R$, on $\fM$, as the generator 
of a one-parameter group of isometries of $\fM$, and so of 
automorphisms of $\Clif(\fM)$, determines a derivation of 
$\Clif(\fM)$, which just extends the action of $R$ on $\fM$.  
Then $L_W^0$ defines a derivation of $\Clif(\cT(G/K))$, which 
we denote again by $L_W^0$.  It is defined by 
\[
(L_W^0\varphi)(x) = L_W^0(x)(\varphi(x)),
\]
where $L_W^0(x)$ here denotes the extension of $L_W^0(x)$ on 
$\fM$ to a derivation of $\Clif(\fM)$.  Since the sum of two 
derivations is a derivation, $\nabla_W^0 + L_W^0$ acts as a 
derivation on $\Clif(\cT(G/K))$.  It is easily checked that 
$\nabla^0 +L^0$ is then an $A$-linear map from $\cT(G/K)$ 
into the algebra of derivations of $\Clif(\cT(G/K))$.  This 
is the ``Levi--Civita Clifford connection'' that is used when 
trying to define the Dirac operator for the 
$\lambda$-invariant Riemannian metric on $G/K$.

We note that for any $L$ which is skew-symmetric the above 
construction works for $L$ in place of $L^0$, that is, the 
construction works for any connection on $\cT(G/K)$ that is 
compatible with the Riemannian metric.  We have not yet used 
the torsion $= 0$ condition.


\section{The Hodge--Dirac operator}
\label{sec8}

To define the usual Dirac operator on $G/K$ for its canonical 
Riemannian metric (which depends on the $\Ad$-invariant inner 
product on $\fG$ that has been chosen), one must deal with 
the issues of whether $G/K$ is $\spin$ or $\spin^c$, and with 
the bookkeeping details coming from whether $G/K$ is of even 
or odd dimension.  We will not pursue these aspects in this 
paper (so see \cite{Frd, GVF, Slb, BIL, Schr}).  But one 
always has the generalized Dirac operator that is often 
called the Hodge--Dirac operator.  (See 9.B of \cite{GVF}, 
except that here we work over $\bR$ rather than $\bC$.) 

We work first with smooth functions, before putting on a 
Hilbert-space structure.  We need a representation of 
$\Clif(\cT(G/K))$ to serve as ``spinors''.  We take the 
left-regular representation of $\Clif(\cT(G/K))$ on itself.  
So our Dirac operator will be an operator on $\cS = 
\Clif(\cT(G/K))$.  We recall that we have required the 
sections of $\Clif(\cT(G/K))$ to be smooth.

Let $\nabla$ be any connection on $\cT(G/K)$ compatible with 
the Riemannian metric on $G/K$, and extend it to 
$\Clif(\cT(G/K))$ as done in the previous section.  For 
$\varphi \in \cS$ define $d\varphi$ by $d\varphi(W) = 
\nabla_W\varphi$ for $W \in \cT(G/K)$.  Thus we can view 
$d\varphi$ as an element of $\cS \otimes \cT^*(G/K)$, where 
$\cT^*(G/K)$ denotes the cross-section module of the 
cotangent bundle.  But by means of the Riemannian metric we 
can identify $\cT^*(G/K)$ with $\cT(G/K)$.  When $d\varphi$ 
is viewed as an element of $\cS \otimes \cT(G/K)$, we denote 
it, with some abuse of notation, by $\grad_{\varphi}$.  Let 
$c$ denote the product on the algebra $\cS = 
\Clif(\cT(G/K))$, viewed as a linear map from $\cS \otimes 
\cS$ to $\cS$.  We view $\cT(G/K)$ as a subspace of $\cS$, 
and so we view $\cS \otimes \cT(G/K)$ as a subspace of $\cS 
\otimes \cS$.  In this way we view $\grad_{\varphi}$ as an 
element of $\cS \otimes \cS$, to which we can apply $c$.  We 
can then define an operator, $D$, on $\cS$ by
\begin{equation}
\label{eq8.1}
D\varphi = c(\grad_{\varphi})
\end{equation}
for all $\varphi \in \cS$.  When $\nabla$ is the Levi--Civita 
connection, this will be our Hodge--Dirac operator, but we do 
not yet assume that $\nabla$ has torsion $0$.

Let us obtain a more explicit formula for $\grad_{\varphi}$, 
and so for $D$.  Consider a standard module frame $\{W_j\}$ 
for $\cT(G/K)$, for example $\{{\hat X}_j\}$ where $\{X_j\}$ 
is an orthonormal basis for $\fG$.  Then for any $V \in 
\cT(G/K)$ we have
\[
d\varphi(V) = \nabla_V\varphi = \nabla_{\sum W_j\langle 
W_j,V\rangle_A} \varphi = \sum(\nabla_{W_j}\varphi)\langle 
W_j,V\rangle_A,
\]
and so we have
\[
\grad_{\varphi} = \sum(\nabla_{W_j}\varphi) \otimes W_j.
\]
When we apply the Clifford product $c$ to this formula for 
$\grad_{\varphi}$, but also use our earlier ``dot'' notation 
for the Clifford product, we find that 
\begin{equation}
\label{eq8.2}
D\varphi = \sum(\nabla_{W_j}\varphi) \cdot W_j.
\end{equation}
We emphasize that $D$ is independent of the choice of 
standard module frame, as can easily be seen directly form the fact that
$(U, V) \mapsto (\nabla_U \varphi) \cdot V$ is clearly $A$-bilinear.
Thus we can choose different frames at 
our convenience to do computations.  A first instance of this 
appears in the next paragraph.

Suppose that the connection $\nabla$ is $\lambda$-invariant.  
We saw in the previous section that this implies that its 
extension to $\cS$ satisfies $\lambda_y(\nabla_W\varphi) = 
\nabla_{\lambda_yW}(\lambda_y\varphi)$.  Let us see that this 
implies that $D$ commutes with $\lambda$.  For a given 
standard module frame $\{W_j\}$ we have
\begin{eqnarray*}
D(\lambda_y\varphi) &= &\sum \nabla_{W_j}(\lambda_y\varphi) 
\cdot W_j = 
\sum(\lambda_y(\nabla_{\lambda_y^{-1}W_j}\varphi)) \cdot 
\lambda_y(\lambda_y^{-1} W_j) \\
&= 
&\lambda_y\left(\sum(\nabla_{\lambda_y^{-1}W_j}\varphi)\right
) \cdot (\lambda_y^{-1}W_j) = \lambda_y(D\varphi),
\end{eqnarray*}
where we have used the easily verified fact that 
$\{\lambda_y^{-1}W_j\}$ is again a standard module frame.  
Thus
\[
D\lambda_y = \lambda_y D
\]
for all $y \in G$.

The elements of $A = C^{\infty}(G/K)$ act as pointwise 
multiplication operators on $\cS$, and it is important to 
calculate their commutators with $D$.  For $f \in A$ let 
$M_f$ denote the corresponding operator on $\cS$.  Then by 
the Leibniz rule for $\nabla$
\begin{eqnarray*}
[D,M_f](\varphi) &= &\sum(\nabla_{W_j}(\varphi f)) \cdot W_j 
- (\sum(\nabla_{W_j}\varphi) \cdot W_j)f \\
&= &\varphi \cdot \left( \sum W_j(\delta_{W_j}f)\right).
\end{eqnarray*}
But much as above
\begin{eqnarray*}
df(V) &= &\delta_Vf = \delta_{\sum W_j\langle 
W_j,V\rangle_A}f \\
&= &\sum (\delta_{W_j}f) \langle W_j,V\rangle_A = 
\left\langle \sum W_j(\delta_{W_j}f),V\right\rangle_A,
\end{eqnarray*}
so that $\sum W_j(\delta_{W_j}f)$ is the usual gradient, 
$\grad_f$, of $f$.  That is:

\setcounter{theorem}{2}
\begin{proposition}
\label{prop8.3}
For $f \in A$ and $\varphi \in \cS$ we have
\[
[D,M_f]\varphi = \varphi \cdot \grad_f,
\]
the product being that in $\Clif(\cT(G/K)) = \cS$.
\end{proposition}

We now want a Hilbert space structure on $\cS$.  As on any 
Clifford algebra, there is a canonical normalized trace, 
$\tau$, on $\Clif(\fM)$, determined by the properties that 
$\tau(1) = 1$ and that if $\{Y_j\}_{j=1}^q$ is a  set of 
mutually orthogonal elements of $\fM$ then $\tau(Y_1\cdot 
{\dots} \cdot Y_q) = 0$.  Note in particular that  $\tau(X 
\cdot Y) = -\langle X,Y\rangle_{\fM}$ for any $X,Y \in \fM$.  
(One way to see the existence of $\tau$ is to consider an 
orthonormal basis for $\fM$ and to decree that $\tau$ is $0$ 
on all of the corresponding basis elements \cite{Frd, GVF} 
for $\Clif(\fM)$ except $1$.)  On $\Clif(\fM)$ there is the 
standard involutory automorphism carrying $X$ to $-X$ for $X 
\in \fM$, and the standard involutory anti-automorphism that 
reverses the order of products.  We let ${}*$ denote the 
composition of these two, so that ${}^*$ is an 
anti-automorphism of $\Clif(\fM)$ that carries $X$ to $-X$ 
for $X \in \fM$.  It is involutory since its square is an 
automorphism that is the identity on $\fM \subset 
\Clif(\fM)$.  On $\Clif(\fM)$ we define an ordinary inner 
product, $\langle \cdot, \ \cdot \rangle_c$, by 
\[
\langle\varphi, \ \psi\rangle_c  = \tau(\varphi^*\cdot\psi)
\]
for $\varphi,\psi \in \Clif(\fM)$.  It is easy to verify that 
for any orthonormal basis for $\fM$ the corresponding basis 
for $\Clif(\fM)$ is orthonormal for this inner product, and 
that the inner product is, in fact, positive definite.  Both 
the left and right regular representations of $\Clif(\fM)$ on 
itself are easily seen to be $*$-representations for this 
inner product.  In particular, elements of $\fM$ act, on left 
and right, as skew-symmetric operators.

We apply all of the above structures to $\cS = 
\Clif(\cT(G/K))$.  We obtain an involution on $\cS$ defined 
by $(\varphi^*)(x) = (\varphi(x))^*$, and we obtain an 
$A$-valued inner product, $\langle\cdot,\cdot\rangle_A$, on 
$\cS$ defined by
\[
\langle\varphi, \ \psi\rangle_A(x) = 
\langle\varphi(x), \ \psi(x)\rangle_c
\]
for $\varphi,\psi \in \cS$.  The left and right regular 
representations of $\Clif(\cT(G/K))$ on $\cS$ are then 
``$*$-representations'', that is, for any $\theta \in \cS$ we 
have
\[
\langle\theta \cdot \varphi, \ \psi\rangle_A = 
\langle\varphi, \ \theta^* \cdot \psi\rangle_A,
\]
and similarly for $\theta$ acting on the right.  Finally, we 
can define an ordinary inner product on $\cS$ by
\[
\langle\varphi, \ \psi\rangle = \int_{G/K} 
\langle\varphi, \ \psi\rangle_A(x)dx.
\]
When $\cS$ is completed for this inner product we obtain our 
Hilbert space of ``spinors'' for the corresponding 
``Hodge-Dirac'' operator.  We denote this Hilbert space by 
$L^2(\cS,\tau)$. We can now view the Hodge-Dirac operator as an
unbounded operator on $L^2(\cS,\tau)$ with domain $\cS$.  

One reason for the importance of the 
torsion $= 0$ condition, or at least a weak version of it, is in determining
whether $D$ is formally self-adjoint, that is,
\[
\langle D\varphi, \ \psi\rangle = \langle\varphi, \ D\psi\rangle \ 
\]
for all $\varphi, \psi \in \cS$.
For any $U \in \cT(G/K)$ let $T_\nabla^U$ be the $A$-endomorphism of
$\cT(G/K)$ defined by $T_\nabla^U(V) = T_\nabla(U, \ V)$. Then $T_\nabla^U$ 
is given by a function on $G$ who values are operators on $\fM$, according
to Proposition \ref{prop3.1}. Thus we can define 
$\mathrm{trace}(T^U_\nabla)$ pointwise as a function in $A$. Equivalently,
$\mathrm{trace}(T^U_\nabla) = \sum_j \langle T_\nabla(U, W_j), W_j\rangle_A$
for one (hence every)
standard module frame $\{W_j\}$ for $\cT(G/K)$.

\begin{theorem}
\label{th8.4}
Let $\nabla$ be any $G$-invariant connection on $\cT(G/K)$ compatible
with our chosen  
Riemannian metric on $G/K$.  
Let $D$ be the Hodge-Dirac operator defined 
as above for $\nabla$, viewed 
as an unbounded operator on $L^2(\cS,\tau)$ with domain 
$\cS$.  Then $D$ is formally self-adjoint if and only if  
\[ 
\mathrm{trace}(T^U_\nabla)  \ \ = \ \ 0
\]
for all $U \in \cT(G/K)$, where $T_\nabla$ is the torsion of $\nabla$.

\end{theorem}

\begin{proof}
As one might suspect, the proof is a complicated version of 
``integration by 
parts'' or the divergence theorem.  We somewhat follow the pattern of the proof
in section 3.2 of \cite{Frd} or of 
proposition~$9.13$ 
of \cite{GVF}.  Let $\{W_j\}$ be a standard module frame for 
$\cT(G/K)$.  
We use first the Leibniz rule for $\nabla$ extended to $\cS$, 
and then the 
fact that $\nabla$ is compatible with the Riemannian metric, 
in order to 
calculate that for $\varphi,\psi \in \cS$ we have
{\setlength\arraycolsep{2pt}   
\begin{eqnarray*}
\langle D\varphi, \ \psi\rangle_A &- 
&\langle\varphi, \ D\psi\rangle_A 
 \ =  \ \sum_j (\langle(\nabla_{W_j}\varphi) \cdot W_j, \ 
\psi\rangle_A - \langle\varphi, \ (\nabla_{W_j}\psi) \cdot 
W_j\rangle_A \\
&= &\sum_j (-\langle\nabla_{W_j}\varphi, \ \psi \cdot 
W_j\rangle_A - \langle\varphi, \ \nabla_{W_j}(\psi\cdot W_j) 
- \psi \cdot (\nabla_{W_j}W_j)\rangle_A \\
&= &\sum_j (-\delta_{W_j}(\langle\varphi, \ \psi \cdot 
W_j\rangle_A) + \langle \varphi, \ \psi \cdot 
(\nabla_{W_j}W_j)\rangle_A).
\end{eqnarray*}
}
For given $\varphi$ and $\psi$ the function $V \mapsto \langle \varphi, \ \psi \cdot V\rangle_A$
is $A$-linear, and so by the self-duality of $\cT(G/K)$ for its Riemannian metric there is a
$U \in \cT(G/K)$ such that $\langle \varphi, \psi \cdot V\rangle_A = \langle U, V\rangle_A$
for all $V \in \cT(G/K)$. By letting $\psi = 1$ and $\varphi = U$  
(viewed as elements in $\Clif(\cT(G/K))$) we see
that any $U$ arises in this way. The above displayed expression is then equal to

\[
\sum_j -\delta_{W_j}(\langle U, W_j\rangle_A) + \langle U, \nabla_{W_j}W_j \rangle_A \ = \ 
-\sum_j \langle \nabla_{W_j}U, W_j\rangle_A \ .
\]
We thus see that $D$ is formally self-adjoint exactly if 
$\int \sum_j \langle \nabla_{W_j}U, W_j \rangle_A \ = \ 0$
for all $U \in \cT(G/K)$ and one, hence every, standard module frame $\{W_j\}$.
 
Now by the definition of $T_\nabla$
\[
\langle \nabla_{W_j} U, W_j\rangle_A \ = 
\ \langle \nabla_U W_j \ - \ T_\nabla (U, W_j) \ - \ [U, W_j] \ , \ W_j\rangle _A   .
\]
Notice then that, by the $A$-bilinearity of the 
inner product, $\sum \langle W_j, W_j \rangle_A$  is
independent of the choice of standard module frame, and that 
$\{W_j(x)\}$ is a frame for $\fM$ for any $x \in G$, so that we can evaluate
the sum by using a frame for $\fM$ that consists of an orthonormal basis
for $\fM$. 
From this we see that 
$\sum \langle W_j, W_j \rangle_A \ \equiv \ \dim(\fM)$.
Consequently, by the compatibility of $\nabla$ with the Riemannian metric,
for any $U \in \cT(G/K)$ we have

\[
0 = \delta_U(\sum_j \langle W_j, W_j\rangle_A) = 
\sum_j \langle \nabla_U W_j, W_j \rangle_A + \langle W_j, \nabla_U W_j\rangle_A =
2\sum_j \langle \nabla_U W_j, W_j\rangle_A      .
\]
Consequently $\sum_j \langle \nabla_U W_j, W_j \rangle_A  =  0$   . (We remark that
this fact depends on the pointwise argument just above, and that the analogous argument can
fail for modules over a non-commutative $A$ that contains proper isometries.) Thus

\[
 \sum_j \langle \nabla_{W_j} U, W_j \rangle_A =
 -\sum_j \langle T_\nabla (U, W_j), W_j \rangle_A \ - \ \sum_j \langle [U, W_j], W_j \rangle_A  \ .
 \]
 Let $\nabla^t = \nabla^0 \ + \ L^0$, the Levi-Civata connection. We can apply the above
 equation to $\nabla^t$ and use that $\nabla^t$ is torsion-free to get an espression for the
 last term above. In this way we find that
 
 \[
 \langle \nabla_{W_j} U, W_j\rangle_A \ = 
 -\sum_j \langle T_\nabla(U, W_j), W_j\rangle_A \ + \ \sum_j\langle \nabla^t_{W_j}U, W_j\rangle_A \ .
 \]
 We will show shortly that
 \[
 \int_{G/K} \sum_j \langle \nabla^t _{W_j}U, W_j \rangle_A \ = \ 0 
 \]
 for all $U$. This will show that $D$ is formally self-adjoint if and only if 
 \[
 0 \ = \ \int \sum_j \langle T_\nabla (U, W_j), W_j \rangle_A \ = \ \int \mathrm{trace}(T^U_\nabla)
 \]
for all $U$. But the
integrand of these latter integrals is clearly $A$-linear in $U$, so when we replace
$U$ by $fU$ for any $f \in A$ the $f$ comes outside the inner product. Since $f$
is arbitrary, this means that the integral is 0 for all $U$ exactly if the integrand itself
is 0 for all $U$, and that is the condition in the statement of the theorem.
 
 Thus we have basically reduced the proof of the theorem to treating $\nabla^t$.
 We remark that, because $\nabla^t$ is the Levi-Civita connection,
 $  \sum_j \langle \nabla^t _{W_j} U, W_j \rangle_A \ = \ 
 \mathrm{div}(U)$, so that in effect we need to prove the divergence
 theorem $\int \mathrm{div}(U) \ = \ 0$.
Since $\nabla^t$ is compatible with the Riemannian metric, we have
\[
\langle \nabla^t_{W_j} U, W_j \rangle_A =
\delta_{W_j}(\langle U, W_j\rangle_A) \ - \ \langle U, \nabla^t_{W_j}W_j \rangle_A \ .
\]
Now for any $X \in \fG$ 
and any $f \in A$ we have $\delta_{\hat X}(f) = \lambda_Xf$, and
\[
\int_G (\lambda_Xf)(x)dx = 0
\]
because $\lambda_Xf$ is the uniform limit of the quotients
\[
(\lambda_{\exp(-tX)}f - f)/t
\]
as $t$ converges to $0$ (see lemma 1.7 of chapter 0 of 
\cite{Tyl}, where 
one needs ${\mathcal L}$ to be invariant), and the integral 
over $G$ of 
each of these quotients is clearly $0$ by the left-invariance 
of Haar measure.  Thus if we 
choose our standard module frame to be $W_j = {\hat X}_j$ for 
an orthonormal 
basis $\{X_j\}$ for $\fG$, we see that for each $j$
\[
\int_{G/K} \delta_{{\hat X}_j}(\langle U,
{\hat X}_j\rangle_A) = 0.
\]
Consequently
\[
\int_{G/K} \sum_j \langle \nabla^t_{{\hat X}_j} U, {\hat X}_j\rangle_A = 
\int_{G/K} \sum_j \langle U, \nabla^t_{{\hat X}_j}{\hat X}_j\rangle_A  \ .
\]
Since we want this to be 0 for all U, it is clear that we need to show that
\[
\sum_j \nabla^t_{{\hat X}_j} {\hat X}_j \ = \ 0 \ .
\]
Now $\nabla^t \ = \ \nabla^0 \ + \ L^0$, and for any $W \in \cT(G/K)$ we have
\[
(L^0_{W}W)(x) \ = \ (1/2)P[W(x), W(x)] \ = \ 0 \ .
\]
Thus we only need to show that $\sum_j \nabla^0_{{\hat X}_j} {\hat X}_j \ = \ 0$.
But
\[
(\sum_j \nabla^0_{{\hat X}_j} {\hat X}_j)(x) \ = \ -\sum_j P[P \Ad^{-1}_x X_j, \ \Ad^{-1}_x X_j] \ .
\]
Note that $\{\Ad^{-1}_x X_j\}$ is an orthonormal basis for $\fG$ for each $x$.
Now $(X, Y) \mapsto P[PX, Y]$ is bilinear, and so it is easily seen that
the value of $\sum_j P[PX_j, X_j]$ does not depend on the choice of
orthonormal basis. We can then choose our basis such that $X_1, \dots X_p$
is a basis for $\fM$ while $X_{p+1}, \dots, X_n$ is a basis for $\fK$. Then
for $j \leq p$ we have $PX_j = X_j$ so that the corresponding terms in the
sum are 0, while for $j > p$ we have $PX_j = 0$ so that again the corresponding terms
in the sum are 0.

\end{proof}

We can use this to show that even though the canonical connection is often not torsion-free,
we have:

\setcounter{theorem}{4}
\begin{corollary}
\label{cor8.5}
The Hodge-Dirac operator for the canonical connection $\nabla^0$ is formally self-adjoint.
\end{corollary}

\begin{proof}

We apply the criterion of Theorem \ref{th8.4} to the canonical connection. 
From our calculation of $T_{\nabla^0}$ in Section 6 we see that for any 
$U \in \cT(G/K)$ we have $T^U_{\nabla^0}(x) \ = \ -P\circ \ad_{U(x)}$ .
We can use an orthonormal basis
$X_1, \dots, X_p$ for $\fM$ in calculating the trace. Extend
this basis by an orthonormal basis $X_{p+1}, \dots, X_n$
for $\fK$. Then for any $U(x) = Y \in \fM$ we have 
\[
\mathrm{trace}(P \circ \ad_Y) = \sum^p \langle P[Y, X_j], X_j \rangle_\fM
=  \sum^n \langle [Y, X_j], X_j \rangle_\fG  = 
\mathrm{trace}(\ad_Y \mathrm{on}  \ \fG )
\]
since $[Y, X_j] \in \fM$ if $X_j \in \fK$ so they are orthogonal. 
But for each $j$ we have $\langle [Y, X_j], X_j\rangle = \langle Y, [X_j, X_j]\rangle = 0$,
and so 
$\mathrm{trace}(\ad_Y \mathrm{on}  \ \fG )  = \ 0$. Thus $\nabla^0$ satisfies
the criterion of Theorem \ref{th8.4}, and so its Hodge-Dirac operator
is formally self-adjoint.
\end{proof}

I have not seen mentioned in the literature the possibility
that some non-torsion-free connections can nevertheless have 
formally self-adjoint Dirac operators.

A simple further calculation using the fact that $L_U$ is a skew-adjoint
operator and so has trace 0, gives:
\begin{corollary}
\label{cor8.6}
Let $\nabla$ be a connection compatible with our chosen Riemannian metric,
and let $L \ = \ \nabla \ - \ \nabla^0$. Then the Hodge-Dirac operator for
$\nabla$ is formally self-adjoint if and only if for one (hence every)
standard module frame $\{W_j\}$ for $\cT(G/K)$ we have
\[ 
\sum_j L_{W_j}W_j\ \ = \ \ 0 \ .
\]
\end{corollary}

For the essential self-adjointness of Dirac operators see, 
for example, 
section~$9.4$ of \cite{GVF} and section~$4.1$ of \cite{Frd}.

Let us now consider the operator norm of $[D,M_f]$ for $f \in 
A$.  We saw earlier that $[D,M_f]\varphi = \varphi \cdot 
\grad_f$ for $\varphi \in \cS$, and that $\grad_f \in \cT(G/K)$.  
Let $R$ denote the right-regular representation of $\cS$ on 
itself and so on $L^2(\cS,\tau)$.  Now the norm of any 
bounded operator, $T$, on an inner-product space satisfies 
the $C^*$-condition $\|T\|^2 = \|T^*T\|$, so for any $\psi 
\in \cS$ we have $\|R_{\psi}\|^2 = \|R_{\psi}^*R_{\psi}\|$.  
When we use this for $\psi = V \in \cT(G/K)$, and recall that 
elements of $\cT(G/K)$ act as skew-symmetric operators for 
the Clifford product, we see that
\[
\|R_V\|^2 = \|R_V^*R_V\| = \|-R_{V\cdot V}\| = \|R_{\langle 
V,V\rangle_A}\|.
\]
But simple arguments show that for any $g \in A$ we have 
$\|R_g\| = \|g\|_{\infty}$.  When we apply all of this for $V 
= \grad_f$ we obtain
\[
\|R_{\grad_f}\|^2 = \|\langle 
\grad_f,\grad_f\rangle_A\|_{\infty},
\]
so that $\|R_{\grad_f}\| = \|\grad_f\|_{\infty}$ for the 
evident meaning of this last term.  But a standard argument 
(e.g., following definition~$9.13$ of \cite{GVF}) shows that 
if we denote by $\rho$ the ordinary metric on a Riemannian 
manifold $M$ coming from its Riemannian metric, then for any 
two points $p$ and $q$ of $M$ we have
\[
\rho(p,q) = \sup\{|f(p) - f(q)|: \|\grad_f\|_{\infty} \le 
1\}.
\]
On applying this to $G/K$, using what we found above for the 
Dirac operator, we obtain, for $\rho$ now the ordinary metric 
on $G/K$,
\[
\rho(p,q) = \sup\{|f(p) - f(q)|: \|[D,M_f]\| \le 1\}.
\]
This is the formula on which Connes focused for general 
Riemannian manifolds \cite{Cn7, Cn3} as it shows that the 
Dirac operator contains all the metric information (and, in 
fact, much more) for the manifold.  This is his motivation 
for advocating that metric data for ``non-commutative 
spaces'' be encoded by providing them with a ``Dirac 
operator''.

We showed earlier that our $D$ commutes with the action 
$\lambda$ of $G$. This is exactly the manifestation in terms 
of $D$ of the fact that the action of $G$ on $G/K$ is by 
isometries for the Riemannian metric and its ordinary metric.

It would be interesting to understand how all of the above relates
to Connes' action principle for selecting \emph{the} Dirac operator from
among all of the spectral triples that give a specified Riemannnian
\cite{Cn3}. (See also theorem 11.2 and section 11.4 of \cite{GVF}.)
Of course, on the face of it Connes' action principle is just for 
spin-manifolds while many homogeneous spaces are not spin.





\def\dbar{\leavevmode\hbox to 0pt{\hskip.2ex 
\accent"16\hss}d}
\providecommand{\bysame}{\leavevmode\hbox 
to3em{\hrulefill}\thinspace}
\providecommand{\MR}{\relax\ifhmode\unskip\space\fi MR } 
\providecommand{\MRhref}[2]{%
\href{http://www.ams.org/mathscinet-getitem?mr=#1}{#2} 
}
\providecommand{\href}[2]{#2}

\end{document}